\documentclass[a4]{article}

\usepackage{amsfonts}
\usepackage{bm}
\usepackage{bbm}
\usepackage{graphicx}

\newcommand{\btheta}{\boldsymbol{\theta}}
\newcommand{\btau}{\boldsymbol{\tau}}

\newcommand{\be}{\boldsymbol{e}}
\newcommand{\bn}{\boldsymbol{n}}
\newcommand{\bq}{\boldsymbol{q}}

\newcommand{\bI}{\boldsymbol{I}}
\newcommand{\by}{\boldsymbol{y}}
\newcommand{\bY}{\boldsymbol{Y}}

\newcommand{\R}{{\mathbb R}}
\newcommand{\Z}{{\mathbb Z}}

\begin{document}

\title{Approximating time to extinction for endemic infection models}
\author{Damian Clancy\\
Department of Actuarial Mathematics and Statistics\\Maxwell Institute for Mathematical Sciences\\Heriot-Watt University\\ Edinburgh\\ EH14 4AS\\ UK\\
d.clancy@hw.ac.uk\\\\
Elliott Tjia\\
University of Liverpool}
\maketitle

\begin{abstract}
Approximating the time to extinction of infection is an important problem in infection modelling.
A variety of different approaches have been proposed in the literature.
We study the performance of a number of such methods, and characterize their performance in terms of simplicity, accuracy, and generality.
To this end, we consider first the classic stochastic susceptible-infected-susceptible (SIS) model, and then a multi-dimensional generalization of this which allows for Erlang distributed infectious periods.
We find that (i)~for a below-threshold infection initiated by a small number of infected individuals, approximation via a linear branching process works well; (ii)~for an above-threshold infection initiated at endemic equilibrium, methods from Hamiltonian statistical mechanics yield correct asymptotic behaviour as population size becomes large; (iii)~the widely-used Ornstein-Uhlenbeck diffusion approximation gives a very poor approximation, but may retain some value for qualitative comparisons in certain cases; (iv)~a more detailed diffusion approximation can give good numerical approximation in certain circumstances, but does not provide correct large population asymptotic behaviour, and cannot be relied upon without some form of external validation (eg simulation studies).
\vspace*{4mm}

\end{abstract}

\section{Introduction}
In modelling endemic infections, a quantity of particular interest is the persistence time until infection dies out of the population.
For the simplest models, it is possible to compute the distribution of this random variable exactly, based on general Markov chain theory.
In particular, the expected persistence time may be expressed as the solution of a linear system of equations.
However, for large populations and for more complicated models, numerical computation of  this exact solution can become very time-consuming, and may also suffer from numerical instability.
Further, it is not straightforward to use the exact solution to address qualitative issues such as, for instance, the effect of greater variability in individual infectious periods upon the expected persistence time in the population.
Hence it remains of great interest to find simple and accurate approximations.

Our aim is to assess a number of approximation techniques in terms of simplicity, accuracy, and generality.
In order to do this, we focus upon relatively simple models, for which it is possible to compute exact values for comparison.
Thus in section~2 we begin by studying the classic susceptible-infected-susceptible (SIS) model of Weiss and Dishon~\cite{WD71}.
Whilst a useful simple model to start from, this classic SIS model is not an entirely typical infection model.
In particular, it possesses a natural 1-dimensional structure (only the number of infected individuals must be kept track of), whereas most infection models are multi-dimensional (eg one may need to keep track of susceptible, infected, and immune individuals, giving a 3-dimensional model).
In section~3 we therefore extend the SIS model to allow for Erlang-distributed infectious periods, generalising the exponentially distributed infectious periods assumed within the classic SIS~model, and resulting in a model which is naturally multi-dimensional, but nevertheless still relatively simple.
The models we study have the further advantage that in each case, rigorous asymptotic results are known giving the expected persistence time as population size becomes large.
Results for the classic SIS model appear in~\cite{AD98}, and the extension to the SIS model with general infections period distribution in~\cite{BBN16}.
These asymptotic results are derived by techniques which seem unlikely to generalize to broader classes of epidemic models, but provide a basis against which to compare the more general approximation methods that we are chiefly interested in. 

We consider four general approximation techniques.
First, we consider the below threshold case, when the basic reproduction number $R_0$ (the average number of new infections caused by a typical infected individual in
an otherwise susceptible population) is less than one.
Then epidemic threshold theory~\cite{B14}, applicable to a wide class of models, tells us that the infectives process may be approximated by a linear branching process, and consequently the mean persistence time of the infection process may be approximated by the mean time to extinction of the branching process.

An alternative, applicable to both the below threshold ($R_0 \le 1$) and above threshold ($R_0 > 1$) cases, is the approach often referred to as `the diffusion approximation'.
Results such as those contained in section~11.3 of~\cite{EK86} justify approximating the infection model (a Markov jump process) by a diffusion process over finite time intervals, in the limit as population size becomes large.
However, in this large population limit, the mean persistence time generally tends to infinity, so that convergence over finite time intervals is not sufficient (see section~11.4 of~\cite{EK86}).
In spite of this lack of rigorous justification, a number of recent papers have nevertheless used the diffusion approximation to approximate mean persistence times for epidemic models, eg~\cite{DV16,WGPAI14}.

In the case $R_0 > 1$, a simpler diffusion approximation is possible, in which the  diffusion approximation referred to above is linearised about the equilibrium point of the deterministic approximation, leading to an Ornstein-Uhlenbeck diffusion process.
This approximation is considerably simpler than the diffusion approximation described previously, and hence has seen widespread use, eg~\cite{AB00,BN10,C05,CP13,LB07,N99,N02}.

Finally, a rather different approach valid for $R_0 > 1$ uses techniques from Hamiltonian statistical mechanics, and has seen much recent attention in the theoretical physics literature~\cite{AM10,DMRH94,EK04,KM08}.

\section{The classic susceptible-infected-susceptible (SIS) model}
\subsection{Approximation techniques}
We consider first the classic stochastic SIS model introduced by Weiss and Dishon~\cite{WD71}.
In this model, infection spreads between members of a closed population of size~$N$ consisting at time $t$ of $S(t)$ susceptible and $I(t)$ infected individuals.
It is sufficient to focus on $I(t)$, since $S(t) = N - I(t)$ for all $t \ge 0$.
The process $\left\{ I(t) : t \ge 0 \right\}$ is a continuous-time Markov chain, with rates of transition given in table~\ref{SISrates}, where $\beta , \gamma > 0$ are the infection and recovery rate parameters, respectively.
\begin{table}
\begin{tabular}{lll} \hline
Event & State transition & Transition rate \\ \hline
Infection of susceptible & $I \to I+1$ & $(\beta/N) I (N-I)$ \\
Recovery of infected & $I \to I-1$ & $\gamma I$ \\ \hline
\end{tabular}
\caption{Transition rates for the classic SIS model.}\label{SISrates}
\end{table}
The basic reproduction number for this model is given by $R_0 = \beta / \gamma$.
The state space consists of an absorbing state at $I=0$ together with the transient communicating class $C = \left\{ 1,2,\ldots,N \right\}$.
Consequently, absorption at state~0 will occur within finite time with probability~1, and our interest is in the random variable
\begin{eqnarray*}
T &=& \inf \left\{ t \ge 0 : I(t) = 0 \right\} .
\end{eqnarray*}
We define
\begin{eqnarray*}
\tau_i &=& E \left[ T \left| I(0) = i \right. \right]
\end{eqnarray*}
to be the expected time to extinction starting from $i$ infected individuals, and write $\btau = \left( \tau_1 , \tau_2 , \ldots , \tau_N \right)$.
Then denoting by $Q$ the transition rate matrix of the process, and by $Q_C$ the rate matrix restricted to the transient states, we have (Norris, 1997~\cite{N97}, theorem 3.3.3) that 
\begin{eqnarray}
Q_C \btau &=& - {\bm 1} . \label{exact}
\end{eqnarray}
Equation~(\ref{exact}) may be solved numerically for $\btau$. 
Alternatively, equation~(5.28) of Norden~(1982)~\cite{N82} gives the explicit solution
\begin{eqnarray}
\tau_i &=& {1 \over \gamma} \left( 
\sum_{m=1}^i \sum_{j=m}^N
{1 \over j}
\left( {R_0 \over N} \right)^{j-m}
{(N-m)! \over (N-j)!}
\right) . \label{Norden}
\end{eqnarray}

When the infection process is above threshold ($R_0 > 1$) then rather than specifying a fixed initial state~$I(0)$, it is arguably of more interest to consider as initial condition a population in which infection has already settled to an endemic level.
Since $\left\{ I(t) \right\}$ is a process on a finite state-space with all non-absorbing states communicating, 
then (Darroch and Seneta, 1967~\cite{DS67})
there exists a unique quasi-stationary distribution $\bq = \left( q_1 , q_2 , \ldots , q_N \right)$ such that, whatever the initial state of the process,
\begin{eqnarray*}
q_i &=& \lim_{t \to \infty} \Pr \left( I(t) = i \left| I(t) > 0 \right. \right)
\mbox{ for } i=1,2,\ldots,N .
\end{eqnarray*}
The quasi-stationary distribution $\bq$ may be found as the unique solution of
\begin{eqnarray}
\bq Q_C &=& - \lambda \bq \mbox{ with } q_1 + q_2 + \cdots + q_N = 1 , \label{QSD}
\end{eqnarray}
where $-\lambda$ is the eigenvalue of $Q_C$ with largest real part.
The time to extinction from quasi-stationarity is then exponentially distributed with mean $\tau_{\bq} = 1/\lambda = 1 / \gamma q_1$.
Note that the quasi-stationary distribution $\bq$ and the mean extinction time~$\tau_{\bq}$ may both be evaluated for any $R_0 > 0$.
However, for $R_0 < 1$ the process is likely to die out before settling to a quasi-stationary state, so that $\bq$ and $\tau_{\bq}$ are of little practical interest.
In contrast, for $R_0 > 1$ the process is likely to settle around the quasi-stationary distribution for a considerable time prior to eventual extinction. 
Hence for $R_0 > 1$ it is natural to consider a population in which infection is endemic, and to study the time to extinction from quasi-stationarity.

The simplest possible approximation to the above model is the deterministic SIS model, defined as follows.
For sufficiently large population size~$N$, the scaled number of infectives $I(t)/N$ may be approximated by the deterministic process $y(t)$ satisfying
\begin{eqnarray}
{dy \over dt} &=& \beta y (1-y) - \gamma y . \label{deterministic_ODE}
\end{eqnarray}
The deterministic system~(\ref{deterministic_ODE}) has equilibria at $y=0$ and $y^* = 1 - \left( 1/R_0 \right)$.
For $R_0 < 1$, the former equilibrium (of extinction) is stable, and $y(t) \to 0$ as $t \to \infty$ for any initial state $y(0) \in [0,1]$.
For $R_0 > 1$, the endemic equilibrium $y^*$ is stable, with $y(t) \to y^*$ as $t \to \infty$ for any initial state $y(0) \in (0,1]$.
In fact, equation~(\ref{deterministic_ODE}) may be solved explicitly: for $\beta \ne \gamma$, starting from $y(0)=y_0$, we have
\begin{eqnarray}
y(t) &=& {(\gamma - \beta) y_0 {\rm e}^{(\beta - \gamma ) t} \over \gamma - \beta + \beta \left( 1 - {\rm e}^{(\beta - \gamma ) t} \right) y_0} . \label{deterministic_solution}
\end{eqnarray}
Consequently, the process does not become extinct within finite time (for $y_0 > 0$), whatever the value of $R_0$.
For $R_0 > 1$, since $y(t) \to y^* > 0$ as $t \to \infty$, the deterministic process cannot provide any useful approximation for mean extinction time.
For $R_0 < 1$, on the other hand, although $y(t) > 0$ for all $t$, it is at least the case that $y(t) \to 0$ as $t \to \infty$.
Noting that we are approximating the discrete, integer-valued process $I(t)$ in terms of the continuous process $y(t)$, it is natural to employ a continuity correction; that is, for any initial state with $y_0 \in [ 1/2N , 1]$ we define
\begin{eqnarray*}
\tau^{\mbox{Det}} \left( N y_0 \right)  &=& \inf \left\{ t \ge 0 : N y(t) \le 0.5 \right\} ,
\end{eqnarray*}
and then for $i=1,2,\ldots,N$ we may suggest $\tau^{\mbox{Det}} ( i )$ as a crude first approximation to $\tau_i$ in the subcritical case $R_0 < 1$.
Solving equation~(\ref{deterministic_solution}) for $t$ with $y(t) = 1/2N$, we find
\begin{eqnarray}
\tau^{\mbox{Det}} \left( y \right) &=& {1 \over \gamma} \left( {1 \over 1 - R_0} \right) \left( \ln \left( 2 y \right) + 
\ln \left( 
{1 - R_0 \left( 1 - {1 \over 2N} \right) \over 1 - R_0 \left( 1 - {y \over N} \right)}
\right) \right) 
\label{Det}
\end{eqnarray}
for $0.5 \le y \le N$.

The rather crude (first order) deterministic approximation may be improved by considering a (second order) diffusion approximation.
That is, the process $I(t)$ may be approximated by the diffusion process $Y(t)$ satisfying
\begin{eqnarray}
dY &=& \left( {\beta \over N} Y (N-Y) - \gamma Y \right) dt + \sqrt{{\beta \over N} Y(N-Y) + \gamma Y} \; dW \label{diffusion}
\end{eqnarray}
where $W(t)$ is a standard Brownian motion.
Since we are again approximating a discrete process with a continuous process, we again apply a continuity correction, and approximate $\tau_i$ by $\tau^{\mbox{Diff}} ( i )$, where for $0.5 \le y \le N$,
\begin{eqnarray*}
\tau^{\mbox{Diff}} \left( y \right) &=& E \left[ \left. \inf \left\{ t \ge 0 : Y(t) \le 0.5 \right\} \right| Y(0) = y \right] .
\end{eqnarray*}

From Gardiner 2009~\cite{G09}, section 5.5.2, we have that $\tau^{\mbox{Diff}} (y)$ satisfies the Kolmogorov backward equation
\begin{eqnarray}
\left( {\beta \over N} y (N-y) - \gamma y \right) {\partial \tau^{\mbox{Diff}} \over \partial y} 
+ {1 \over 2} \left( {\beta \over N} y (N-y) + \gamma y \right) {\partial^2 \tau^{\mbox{Diff}} \over \partial y^2} &=& -1 \label{Diffusion_KBE}
\end{eqnarray}
for $y \in (0.5,N)$, with boundary conditions $\tau^{\mbox{Diff}} (0.5) = 0$, $\displaystyle \left. {\partial \tau^{\mbox{Diff}} \over \partial y} \right|_{y=N} = 0$ (where we have imposed a reflecting boundary at $y=N$).
The solution to equation~(\ref{Diffusion_KBE}) may be written explicitly (Gardiner~2009~\cite{G09}, equation~5.5.24) as
\begin{eqnarray}
\tau^{\mbox{Diff}} (y) &=& 2 N \int_{1/2N}^{y/N} du \int_{u}^1 {\exp \left( 2 N \int_u^z {\beta v(1-v) - \gamma v \over \beta v(1-v) + \gamma v} \, dv \right) \over \beta z(1-z) + \gamma z} \, dz \nonumber \\
&=& {2 N \over \gamma} \int_{1/2N}^{y/N} du \int_{u}^1
{\left( {1 + R_0 ( 1 - z ) \over 1 + R_0 ( 1 - u )} \right)^{4N \over R_0}
{\rm e}^{2N(z-u)}
\over z \left( 1 + R_0 (1-z) \right)}
\, dz .
\label{Diff}
\end{eqnarray}
The above integral formula for $\tau^{\mbox{Diff}} (y)$ may be evaluated numerically using Matlab.
In the case $R_0 > 1$, we will approximate the time to extinction from quasi-stationarity $\tau_{\bm{q}}$ using $\tau^{\mbox{Diff}} \left( \lfloor N y^* \rfloor \right)$, where $\lfloor N y^* \rfloor$ denotes the integer part of the re-scaled deterministic endemic equilibrium point.

In the supercritical case ($R_0 > 1$), when studying fluctuations around endemic quasi-equilibrium, it is standard (see, for instance,~\cite{N99}) to further approximate the diffusion $Y(t)$ in terms of an Ornstein-Uhlenbeck process centred at the deterministic equilibrium point~$N y^*$.
The drift and diffusion coefficients in equation~(\ref{diffusion}) are approximated to leading order around $Y = N y^*$, yielding the Ornstein-Uhlenbeck process $\tilde Y(t)$ satisfying
\begin{eqnarray*}
d \tilde Y &=& - ( \beta - \gamma ) \left( \tilde Y - N y^* \right) dt + \sqrt{2 N \gamma y^*} \; dW .
\end{eqnarray*}
The stationary distribution of an Ornstein-Uhlenbeck process with drift coefficient $J$ and local variance $\tilde B$, centred at $Ny^*$, is Gaussian with mean $N y^*$ and variance $V$ satisfying
\begin{eqnarray*}
2 J V &=& - \tilde B .
\end{eqnarray*}
In this case we have $J = - ( \beta - \gamma )$ and $\tilde B = 2 N \gamma y^*$, so that $V = N/R_0$.
The quasi-stationary distribution $\bq$ may thus be approximated by a Gaussian distribution with mean $N y^*$, variance $N/R_0$.
Recalling that the mean time to extinction is given by $\tau_{\bq} = 1/ \gamma q_1$, we have the approximation
\begin{eqnarray*}
q_1 &\approx& \sqrt{{R_0 \over 2 \pi N}} \, \exp \left( - {R_0 \over 2 N} \left( N y^* - 1 \right)^2 \right) .
\end{eqnarray*}
For large $N$, recalling that $y^* = 1 - (1/R_0)$ and making the slight further approximation $Ny^* - 1 \approx N y^*$ gives
\begin{eqnarray*}
q_1 &\approx& \sqrt{{R_0 \over 2 \pi N}} \, \exp \left( - {\left( R_0 - 1 \right)^2 \over 2 R_0} N \right) .
\end{eqnarray*}

For $N$ large and $R_0 > 1$, the mean time to extinction from quasi-stationarity $\tau_{\bm{q}}$ may thus be approximated by
\begin{eqnarray}
\tau^{\mbox{OU}} &=& {1 \over \gamma} \sqrt{ 2 \pi N \over R_0} \, \exp \left( {\left( R_0  - 1 \right)^2 \over 2 R_0} N \right) . \label{OU}
\end{eqnarray}

The above argument can only be expected to give a reasonable approximation provided the Gaussian approximation to the quasi-stationary distribution assigns negligible probability to negative numbers of infected individuals.
A rule of thumb suggested by N{\aa}sell~\cite{N99} is to require the coefficient of variation of the Gaussian distribution to be at most $1/3$, on the basis that a Gaussian distribution assigns negligible probability to values more than 3~standard deviations from the mean; that is, we require $N > 9 R_0 / \left( R_0 - 1 \right)^2$.

Although a number of authors have employed diffusion and Ornstein-Uhlenbeck approximations in studying extinction time for infection models, such approximations are most appropriate in describing moderate deviations from the mean, whereas extinction arises as a result of a large deviation.
Consequently, a more appropriate technique may be the methodology from Hamiltonian statistical mechanics described in, for example,~\cite{EK04,KM08}.
Defining $M ( \theta ; t ) = E \left[ {\rm e}^{\theta I(t)} \right]$ to be the moment generating function of $I(t)$, then the Kolmogorov forward equations for the process may be written as
\begin{eqnarray}
{dM \over dt} &=& \beta \left( {\rm e}^{\theta} - 1 \right) \left( {\partial M \over \partial \theta} - {1 \over N} {\partial^2 M \over \partial \theta^2} \right) + \gamma \left( {\rm e}^{-\theta} - 1 \right) {\partial M \over \partial \theta} . \label{mgf_PDE}
\end{eqnarray}
In terms of the cumulant generating function $K(\theta;t) = \ln M(\theta;t)$, equation~(\ref{mgf_PDE}) is equivalent to
\begin{eqnarray}
{dK \over dt} &=&
\beta \left( {\rm e}^{\theta} - 1 \right)
\left( {\partial K \over \partial \theta}
- {1 \over N} \left( {\partial ^2 K \over \partial \theta^2}  + \left( {\partial K \over \partial \theta} \right)^2 \right)
\right) + \gamma \left( {\rm e}^{-\theta} - 1 \right) {\partial K \over \partial \theta}
 . \label{cgf_PDE}
\end{eqnarray}
Supposing ({\em ansatz}) that $K ( \theta ; t ) = N S(\theta;t) + o(N)$ for some function $S(\theta;t)$, then equation~(\ref{cgf_PDE}) becomes
\begin{eqnarray*}
{dS \over dt} &=& \beta \left( {\rm e}^{\theta} - 1 \right)
{\partial S \over \partial \theta}
\left( 1
- {\partial S \over \partial \theta} 
\right)
+ \gamma \left( {\rm e}^{-\theta} - 1 \right) {\partial S \over \partial \theta}
+ o(1) .
\end{eqnarray*}
To leading order, this has the form of a Hamilton-Jacobi equation corresponding to the Hamiltonian
\begin{eqnarray}
H (y, \theta ; t) &=& \beta y(1-y) \left( {\rm e}^{\theta} - 1 \right) + \gamma  y \left( {\rm e}^{-\theta} - 1 \right) . \label{Hamiltonian}
\end{eqnarray}
In this formulation, $y$ represents the scaled number of infected individuals, and  $\theta$ represents (in the terminology of Hamiltonian mechanics) a conjugate momentum variable.
Note that $\left. {\partial S \over \partial \theta} \right|_{\theta=0} = E \left[ I/N \right] + o(1)$.

References~\cite{EK04,KM08} formulate the problem in terms of the probability generating function; we find the above formulation in terms of the cumulant generating function to be more natural. 
An alternative approach is to work directly with the Kolmogorov forward equations in terms of the quasi-stationary probabilities $q_i$, see for example~\cite{DMRH94}. 
The same Hamiltonian~(\ref{Hamiltonian}) is thus obtained.

The Hamiltonian describes a deterministic motion in $\left( y , \theta \right)$ space that follows the equations of motion
\begin{eqnarray*}
{dy \over dt} = {\partial H \over \partial \theta} , \hspace{2cm} {d \theta \over dt} = - \,  {\partial H \over \partial y} .
\end{eqnarray*}

For the specific Hamiltonian~(\ref{Hamiltonian}) the equations of motion are therefore
\begin{eqnarray}
\left.
\begin{array}{rcl}
\displaystyle {dy \over dt} &=& \beta y (1-y) {\rm e}^{\theta}- \gamma y {\rm e}^{-\theta} , \\\\
\displaystyle {d \theta \over dt} &=& - \beta \left( 1 - 2y \right) \left( {\rm e}^{\theta} - 1 \right) - \gamma \left( {\rm e}^{-\theta} - 1 \right)  .
\end{array} \right\} \label{Hamiltonian_system}
\end{eqnarray}
Trajectories of the deterministic SIS model~(\ref{deterministic_ODE}) may be obtained by taking $\theta = 0$ in equations~(\ref{Hamiltonian_system}).
The system~(\ref{Hamiltonian_system}) possesses classical equilibrium points at $(y,\theta) = (0,0)$ and $(y,\theta) = \left( y^* , 0 \right)$, together with a non-classical disease-free equilibrium point at $(y,\theta) = \left( 0 , - \ln \left( R_0 \right) \right)$.

The expected time to extinction starting from the deterministic equilibrium point  satisfies (see~\cite{EK04}) $\ln \left( \tau_{\lfloor N y^* \rfloor} \right) / N \to A$ as $N \to \infty$, where $A$ is known as the `action'.
The value of $A$ is found by integrating along the trajectory of~(\ref{Hamiltonian_system}) that goes from initial state~$(y,\theta) = \left( y^* , 0 \right)$ at time $t = -\infty$ to final state~$(y,\theta) = \left( 0 , - \ln \left( R_0 \right) \right)$ at time $t=+\infty$, with
\begin{eqnarray*}
A &=& \int_{-\infty}^\infty \theta \, {dy \over dt} \, dt .
\end{eqnarray*}

Now $H ( y , \theta )$ is a constant of the motion (sometimes referred to as the `energy' of the system), and $H \left( y^* , 0 \right) = 0$, so along this trajetory we have
\begin{eqnarray*}
\gamma y \left( {\rm e}^{-\theta} -1 \right) &=& \beta y (1-y) \left( 1 - {\rm e}^\theta \right) .
\end{eqnarray*}
It follows that $y=0$, or $\theta = 0$, or
$\theta = - \ln \left( R_0 (1-y) \right)$.
Along the relevant trajectory we do not have $y=0$ or $\theta = 0$ except at the endpoints, and so
\begin{eqnarray*}
A &=& \int_{y^*}^0 \theta \, dy 
= \int_0^{y^*} \ln \left( R_0 (1-y) \right) \, dy 
= \left( 1 / R_0 \right) - 1 + \ln R_0 .
\end{eqnarray*}

For $R_0 > 1$ and $N$ large, using $\tau_{\lfloor N y^* \rfloor}$ to approximate $\tau_{\bm{q}}$, we therefore have that $\ln \left( \tau_{\bq} \right) / N$ may be approximated by $\ln \left( \tau_H \right) / N$, where
\begin{eqnarray}
\tau_H &=& \exp \left( N \left( \left( 1 / R_0 \right) - 1 + \ln R_0 \right) \right) .
\label{H}
\end{eqnarray}

The methods outlined above are quite general, and can potentially be applied to any model for transmission of infection (or other population processes).
More precise and rigorous results specific to the classic SIS model are presented in~\cite{KL89,AD98,DSS05}.
The relevant results that we will use for comparison with the general approximation methods are as follows.
\begin{enumerate}
\item[(i)] For $R_0 < 1$ and $I(0) = i$ fixed, then as $N \to \infty$ the mean extinction time 
$\tau_i$ converges to the extinction time of a linear birth-death process with birth and death rate parameters $\beta , \gamma$, respectively~\cite{AD98}.
That is, $\tau_i \to \tau_i^{\mbox{Lin}}$, where from equation~(7.10), chapter~4, of Karlin and Taylor~\cite{KT75}, we find
\begin{eqnarray}
\tau_i^{\mbox{Lin}} &=& {1 \over \gamma} \left( {1 \over 1 - R_0} \right) 
\left(
\left( 1 - R_0^{-i} \right) \ln \left( 1 - R_0 \right) 
+ \sum_{m=1}^{i-1} {1 - R_0^{m-i} \over m}
\right) .
\label{Lin}
\end{eqnarray}
\item[(ii)] For $R_0 < 1$ and $i/N = I(0) / N \to y_0 \in (0,1]$ as $N \to \infty$, then as $N \to \infty$ we have 
$\tau_i - \tau^{\mbox{DSS}} \left( N y_0 \right) = O(1)$, where~\cite{DSS05}
\begin{eqnarray}
\tau^{\mbox{DSS}} \left( y \right) \ &=& {1 \over \gamma} \left( {1 \over 1 - R_0}  \right) \ln y . \label{DSS}
\end{eqnarray}
\item[(iii)] For $R_0 > 1$, the expected time to extinction starting from quasi-stationarity satisfies $\tau_{\bq} \sim \tau^{\mbox{AD}}$, where~\cite{AD98}
\begin{eqnarray}
\tau^{\mbox{AD}} &=& {1 \over \gamma} \sqrt{2\pi \over N} \, {R_0 \over \left( R_0 - 1 \right)^2} \, \exp \left( N \left( \left( 1 / R_0 \right) - 1 + \ln R_0 \right) \right) . \label{AD}
\end{eqnarray}
\end{enumerate}

Result~(i) for the case $R_0 < 1$ is based on the epidemic threshold theorem, and can be generalised to a wide range of infection models.
In contrast, the derivation of result~(iii) for the case $R_0 > 1$ relies upon the 1-dimensional structure of the model, so that while similar results may apply to more general birth-death processes, it does not seem straightforward to apply such methods to more general, multi-dimensional, infection models.
Doering et al.~\cite{DSS05}, by different arguments, confirmed the formula~(\ref{AD}) for the case $R_0 > 1$.
In the case $R_0 < 1$ with $I(0) / N \to y_0 \in (0,1]$ as $N \to \infty$, however, formula~(\ref{DSS}) of~\cite{DSS05} is in disagreement with a corresponding result from~\cite{KL89}, cited in~\cite{AD98}.
For completeness, we include this alternative result, which states that 
$\tau_i - \tau^{\mbox{KL}} \left( N y_0 \right) \to 0$ with
\begin{eqnarray}
\tau^{\mbox{KL}} \left( y \right) &=& {1 \over \gamma} \left( {\ln \left( y \right) + \ln \left( 1 - R_0 \left( 1 - (y/N) \right) \right) + \gamma^E \over 1 - R_0 \left( 1 - (y/N) \right)} \right) , \label{KL}
\end{eqnarray}
where $\gamma^E$ denotes Euler's constant $\gamma^E \approx 0.577216$.

\subsection{Numerical results}
\subsubsection{Below threshold case $R_0 < 1$}
For all numerical work we fix $\gamma = 1$, since varying $\gamma$ for fixed $R_0$ simply amounts to a scaling of the time axis.

\begin{center}
\begin{figure}
\includegraphics[width=12cm]{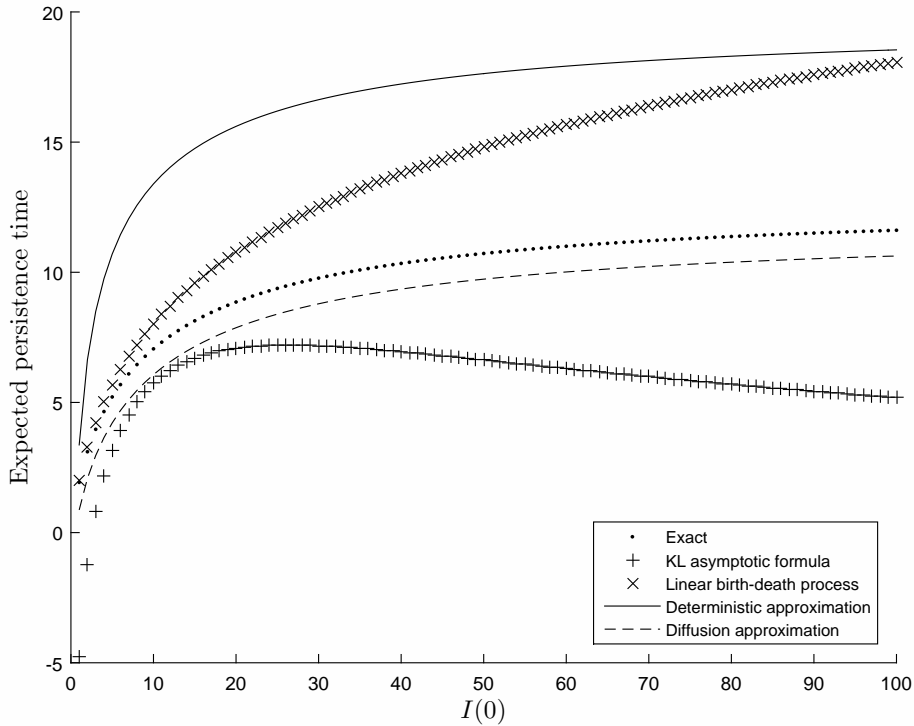}
\caption{Expected persistence time for the classic SIS model from any initial state $I(0)$, together with approximations.
Parameter values $N=100$, $\beta = 0.8$, $\gamma = 1$ (so that $R_0 = 0.8$).
Exact values computed from formula~(\ref{Norden}); KL asymptotic formula~(\ref{KL}); linear birth-death process approximation~(\ref{Lin}); deterministic approximation~(\ref{Det}); diffusion approximation~(\ref{Diff}).}
\label{fig1}
\end{figure}
\end{center}

Figure~\ref{fig1} illustrates how our various approximations compare with the exact mean persistence time $\tau_i$ from any initial state $i \in \left\{ 1,2,\ldots,N \right\}$ for relatively small population size $N=100$, with $R_0 = 0.8$.
We see that the deterministic approximation $\tau^{\mbox{Det}}$ is, as one might expect, least accurate.
The diffusion approximation $\tau^{\mbox{Diff}}$ seems to perform best of all our approximations here.
The approximation $\tau^{\mbox{Lin}}$ performs very well for small initial infective numbers~$i$, but much less well as $i$ becomes large.
Again this is to be expected, since this approximation is based on an assumption that the number of available susceptible individuals remains approximately equal to~$N$ at all times.
The approximation $\tau^{\mbox{KL}}$ performs well in the range $5 \le i \le 20$, but increasingly badly as~$i$ increases, and actually produces negative answers for $i=1,2$.
It is also notable that $\tau^{\mbox{KL}}$, in contrast to the other approximations, is not monotonically increasing in $i$, whereas the exact solution~(\ref{Norden}) clearly is.

\begin{center}
\begin{figure}
\includegraphics[width=12cm]{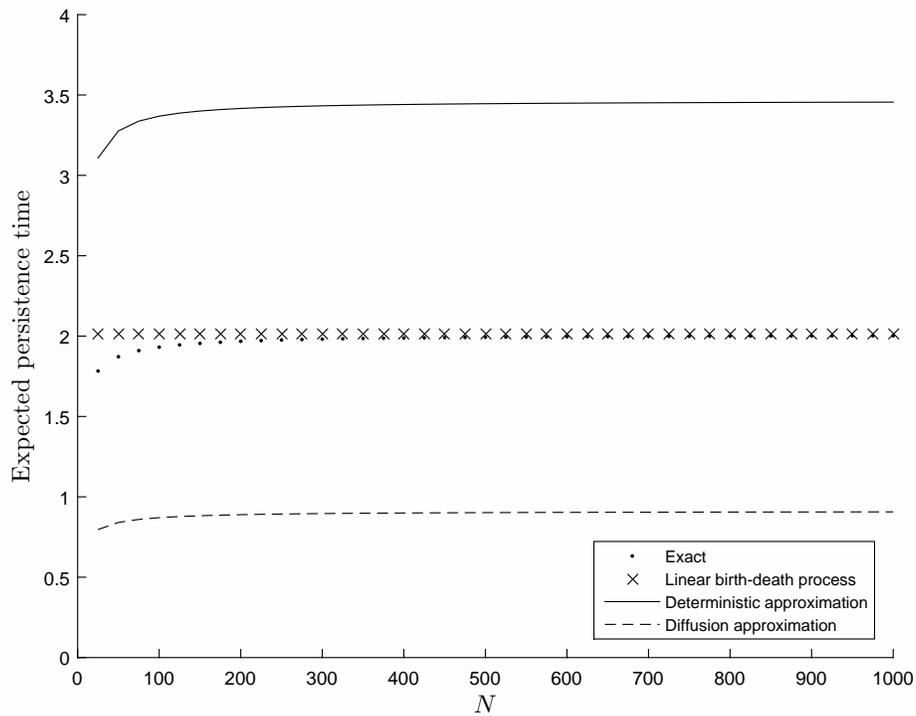}
\caption{Expected persistence time for the classic SIS model starting with $I(0) = 1$ initial infective in a population of size~$N$, together with approximations.
Parameter values $\beta = 0.8$, $\gamma = 1$ (so that $R_0 = 0.8$).
Exact values computed from formula~(\ref{Norden}); linear birth-death process approximation~(\ref{Lin}); deterministic approximation~(\ref{Det}); diffusion approximation~(\ref{Diff}).}
\label{fig2}
\end{figure}
\end{center}

Next, we consider the performance of our approximations across a range of population sizes~$N$, up to the moderately large size $N=1000$.
There are essentially two initial conditions of interest --- either a fixed number of initial infectives, or a fixed fraction of the population initially infected.
Figure~\ref{fig2} shows the performance of relevant approximations for the case of a fixed number of initial infectives, specifically $i=1$.
We see that the approximation $\tau^{\mbox{Lin}}$, whose value is independent of $N$, performs well even for small population sizes~$N$, and very well for larger $N$~values.
Neither $\tau^{\mbox{Det}}$ nor $\tau^{\mbox{Diff}}$ perform nearly as well, and they do not improve as $N$ increases, which is not surprising, since both the deterministic and the diffusion process assume that the number of infective individuals is large enough to be treated as a continuous variable.
For these parameter values, $\tau^{\mbox{KL}}$ takes negative values, and so is not shown.

\begin{center}
\begin{figure}
\includegraphics[width=12cm]{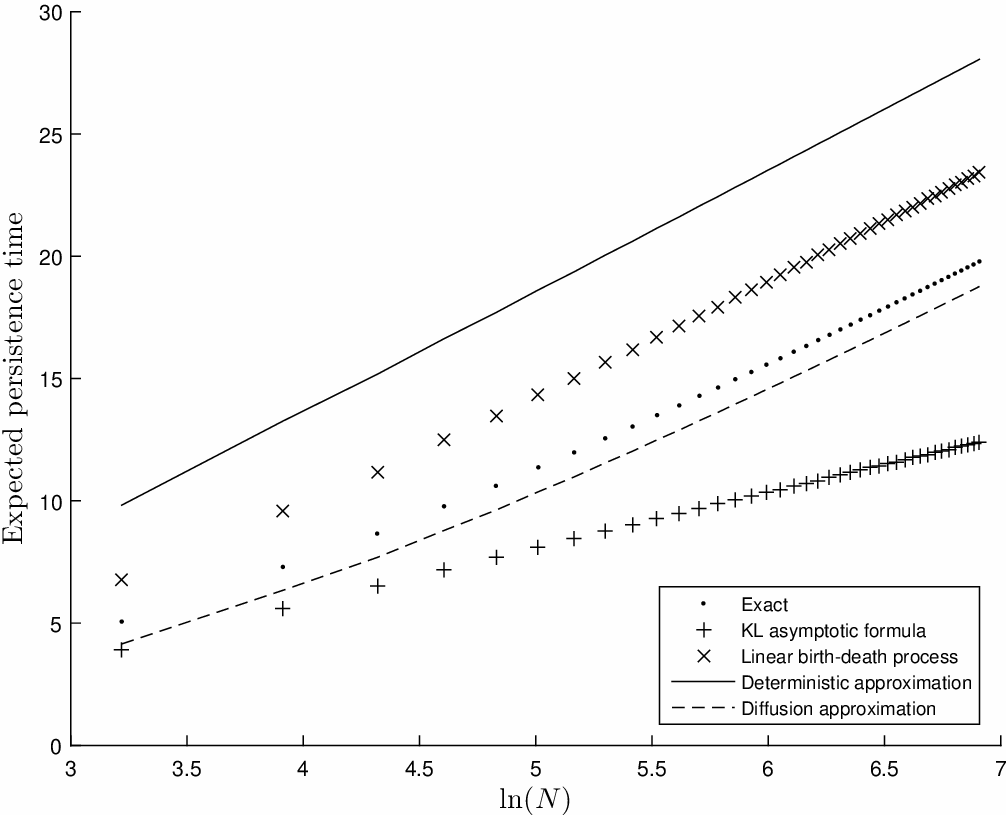}
\caption{Expected persistence time for the classic SIS model plotted against $\ln(N)$ starting with $I(0) = \lfloor 0.3 N \rfloor$ initial infectives in a population of size~$N$, together with approximations.
Parameter values $\beta = 0.8$, $\gamma = 1$ (so that $R_0 = 0.8$).
Exact values computed from formula~(\ref{Norden}); KL asymptotic formula~(\ref{KL}); linear birth-death process approximation~(\ref{Lin}); deterministic approximation~(\ref{Det}); diffusion approximation~(\ref{Diff}).}
\label{fig4}
\end{figure}
\end{center}

For the case of a fixed initially infected fraction $y_0$, figure~\ref{fig4} shows expected persistence time plotted against $\ln (N)$.
The best approximation here is $\tau^{\mbox{Diff}}$, with $\tau^{\mbox{Det}} , \tau^{\mbox{Lin}} , \tau^{\mbox{KL}}$ all performing rather poorly.
In terms of asymptotic (large $N$) behaviour, formula~(\ref{DSS}) of~\cite{DSS05} suggests that, for large~$N$, we should observe a straight line of gradient $(1/\gamma) \left( 1 - R_0 \right)^{-1}$, whereas formula~(\ref{KL}) of~\cite{KL89} predicts a line of gradient $(1/\gamma) \left( 1 - R_0 \left( 1 - y_0 \right) \right)^{-1}$.
From figure~\ref{fig4} it is clear that formula~(\ref{KL}), $\tau^{\mbox{KL}}$, does not predict the correct asymptotic gradient.
The approximations $\tau^{\mbox{Det}}$, $\tau^{\mbox{Diff}}$, $\tau^{\mbox{Lin}}$ all appear to display the correct asymptotic gradient, but with different offsets, with $\tau^{\mbox{Diff}}$ performing best.
Note that the deterministic approximation~(\ref{Det}) predicts the same asymptotic gradient as~\cite{DSS05}.

Overall, our numerical results for the below threshold case suggest that:
(i)~the approximation $\tau^{\mbox{KL}}$ of~\cite{KL89} gives the wrong asymptotic behaviour, and overall does not seem a reliable approximation;
(ii)~the formulae $\tau^{\mbox{DSS}}$ of~\cite{DSS05} and $\tau^{\mbox{Det}}$ give the correct first-order asymptotic behaviour, but $\tau^{\mbox{DSS}}$ differs from the exact answer by an unspecified offset, while $\tau^{\mbox{Det}}$ differs by a substantial constant offset, so that neither provides a useful numerical approximation;
(iii)~in the case of an initial small trace of infection, the best approximation is that derived from an approximating linear birth-death process, $\tau^{\mbox{Lin}}$;
(iv)~in the case of a substantial fraction of the population being initially infected, the best approximation is that derived from an approximating diffusion process, $\tau^{\mbox{Diff}}$.

\subsubsection{Above threshold case $R_0 > 1$}
In the above threshold case, we focus upon the mean persistence time starting from endemicity; that is, we are interested in approximating $\tau_{\bm{q}}$.

Figure~\ref{fig5} shows the values of $\tau_{\bm{q}}$ and relevant approximations for a system not far above threshold ($R_0 = 1.1$).
The best approximation here is provided by the diffusion process, $\tau^{\mbox{Diff}}$.
Note that for the diffusion approximation, we take initial condition $Y(0) = N y^*$, ie starting from the endemic equilibrium point of the deterministic system.
The Ornstein-Uhlenbeck approximation is seen to perform very poorly, and to become even less accurate as $N$ increases.
The approximation $\tau^{\mbox{AD}}$ performs reasonably well for $N \ge 100$, although not as well as $\tau^{\mbox{Diff}}$.

In figure~\ref{fig6} we consider a system well above threshold ($R_0 = 1.5$).
We see that $\tau^{\mbox{AD}}$ gives the correct leading-order asymptotic behaviour; $\tau^{\mbox{Diff}}$ underestimates the mean persistence time; and $\tau^{\mbox{OU}}$ provides a very poor approximation.
Figure~\ref{fig6} reproduces elements of figure~6 of Doering et al.~\cite{DSS05}, although Doering et al.~\cite{DSS05} did not consider the Ornstein-Uhlenbeck approximation, and computed mean time to extinction from a fixed state close to the endemic equilibrium $N y^*$, rather than from quasi-stationarity.
There are also two differences in how the diffusion approximation is computed: firstly, we compute $\tau^{\mbox{Diff}}$ exactly using formula~(\ref{Diff}), whereas Doering et al.~\cite{DSS05} use an asymptotic approximation; secondly, our formula~(\ref{Diff}) employs a continuity correction, stopping the diffusion at $Y = 0.5$, whereas in~\cite{DSS05} the diffusion is simply absorbed at the boundary $Y=0$.
However, these refinements do not make any appreciable difference to the values plotted.

The asymptotic result derived by Doering et al.~\cite{DSS05} for the diffusion approximation (absorbed at $Y=0$) is that, for $R_0 > 1$, the expected persistence time $\tau^{\mbox{Diff}} (y)$ from any initial state $y = O(N)$ satisfies, as $N \to \infty$,
\begin{eqnarray}
\lefteqn{\tau^{\mbox{Diff}} (y)} \nonumber \\
&&\hspace*{-1cm} 
\sim {1 \over \gamma} 
\sqrt{2\pi \over N} 
{R_0 (R_0 + 1) \over 2  ( R_0 - 1 )^2 \sqrt{R_0}}
\exp \left( 
N \left( 
2 \left( {R_0 - 1 \over R_0} \right)
 + {4 \over R_0} 
 \ln \left( {2 \over R_0 + 1} \right)
\right) 
\right) .
\label{FPE}
\end{eqnarray}

As pointed out by Doering et al.~\cite{DSS05}, comparison of formulae~(\ref{FPE}) and~(\ref{AD}) shows that the diffusion approximation does not give the correct asymptotic behaviour of $\ln ( \tau ) / N$, confirming what we see in figure~\ref{fig6}.
To look more closely at the failure of the diffusion approximation, in figure~\ref{fig8} we plot the exponential constants from formulae~(\ref{AD}, \ref{FPE}).
We see that as $R_0$ increases, leading order behaviour of $\tau^{\mbox{Diff}}$ becomes an increasingly poor approximation to the true leading order behaviour of $\tau_{\bm{q}}$.

It is similarly clear from formula~(\ref{OU}) that $\tau^{\mbox{OU}}$ does not give correct leading-order asymptotic behaviour.
On the other hand, the Hamiltonian approach, formula~(\ref{H}), does yield correct leading-order asymptotic behaviour.

In summary, for $R_0$ slightly above~1 and moderate population size~$N$, the best approximation appears to be provided by the diffusion approximation.
However, the diffusion approximation has the wrong leading-order asymptotic behaviour as $N \to \infty$, and gets worse as $R_0$ increases.
Hence this approximation can only be used provided that $N$ is not too large and $R_0$ is not too far above~1.
This is in line with comments of Doering et al.~\cite{DSS05}.
Correct leading-order asymptotic behaviour may be found using the Hamiltonian approach.
However, this approach does not provide a useful numerical approximation, since the result we obtain is that $\tau_{\bm{q}} \sim C ( N , R_0 ) \tau_H$  with $C( N , R_0 ) = o \left( {\rm e}^N \right)$ as $N \to \infty$.
Since the prefactor $C ( N , R_0 )$ is unknown, we cannot use this asymptotic result to evaluate a numerical approximation to $\tau_{\bm{q}}$.
The Ornstein-Uhlenbeck approximation, though widely used, is seen to provide an extremely poor approximation to $\tau_{\bm{q}}$.

\begin{center}
\begin{figure}
\includegraphics[width=12cm]{extinction_times_fig11}
\caption{Expected persistence time for the classic SIS model starting from quasi-stationarity $\left( \tau_{\bm{q}} \right)$ in a population of size $N$, together with approximations.
Parameter values $\beta = 1.1$, $\gamma = 1$ (so that $R_0 = 1.1$ and $y^* \approx 0.0909$).
Exact values computed from formula~(\ref{QSD}); AD asymptotic formula~(\ref{AD}); diffusion approximation~(\ref{Diff}); Ornstein-Uhlenbeck approximation~(\ref{OU}).}
\label{fig5}
\end{figure}
\end{center}

\begin{center}
\begin{figure}
\includegraphics[width=12cm]{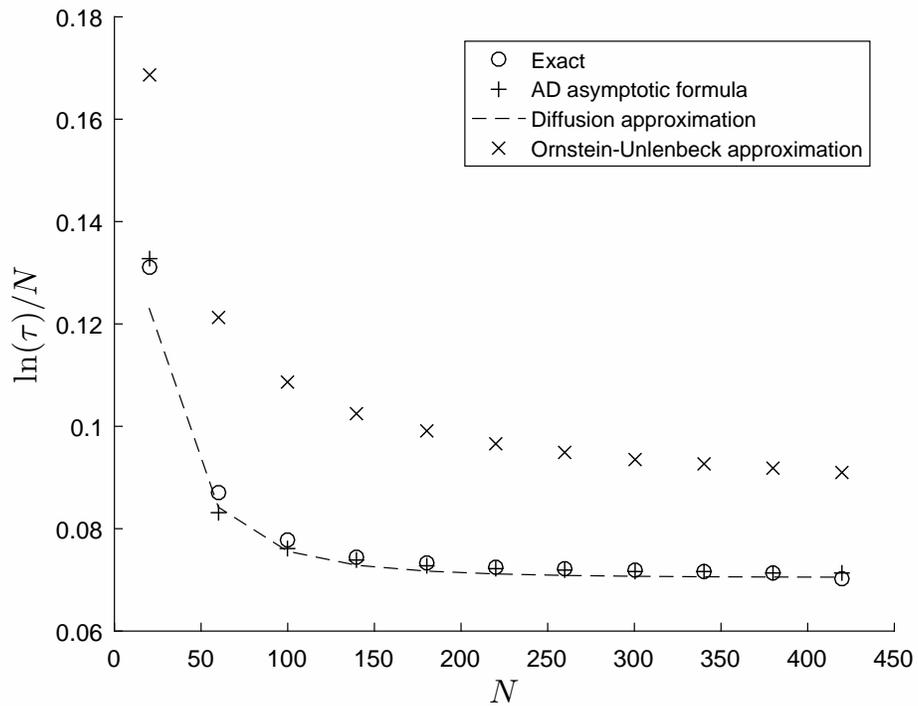}
\caption{$\ln ( \tau_{\bm{q}} ) / N$, where $\tau_{\bm{q}}$ is  the expected persistence time for the classic SIS model starting from quasi-stationarity in a population of size~$N$, together with approximations.
Parameter values $\beta = 1.5$, $\gamma = 1$ (so that $R_0 = 1.5$ and $y^* \approx 0.3333$).
Exact values computed from formula~(\ref{QSD}); AD asymptotic formula~(\ref{AD}); diffusion approximation~(\ref{Diff}); Ornstein-Uhlenbeck approximation~(\ref{OU}).}
\label{fig6}
\end{figure}
\end{center}

\begin{center}
\begin{figure}
\includegraphics[width=12cm]{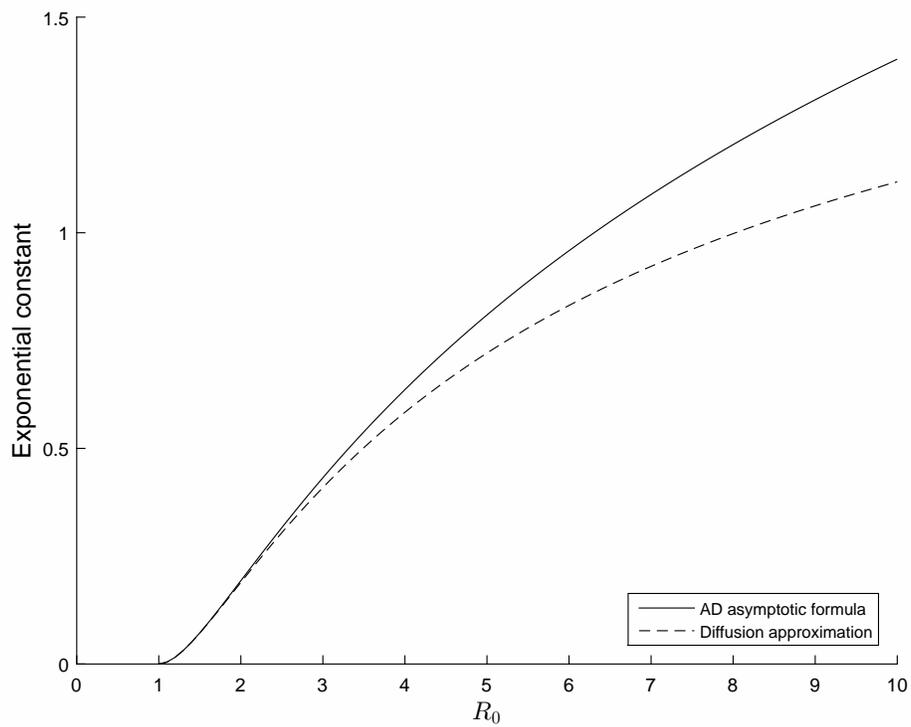}
\caption{Exponential constants appearing in the asymptotic formula $\tau^{\mbox{AD}}$ of~\cite{AD98} for the SIS model and in the asymptotic formula~(\ref{FPE}) of~\cite{DSS05} for the diffusion approximation.}
\label{fig8}
\end{figure}
\end{center}

\section{Generalising the infectious period distribution}
\subsection{Approximation techniques}
The classic SIS model makes the implicit assumption that each individual's infectious period follows an exponential distribution.
This is purely a mathematical convenience, not motivated by biological realism.
We can improve the biological realism of the model by allowing infectious periods to follow an Erlang distribution, using the `method of stages'.
That is, when an individual becomes infected, it passes through $k$ infectious stages, remaining in each stage for an exponentially distributed time of mean $1 / k \gamma$, before returning to susceptibility.
Thus the total infectious period follows an Erlang distribution, and the process can be modelled as a Markov chain $\left\{ \left( I_1 (t) , I_2 (t) , \ldots , I_k (t) \right) : t \ge 0 \right\}$, where $I_m (t)$ is the number of individuals in infectious stage~$m$ at time~$t$, and the number of susceptible individuals at time $t$ is $S(t) = N - \sum_{m=1}^k I_m (t)$.
Transition rates are given in table~\ref{SIS_Erlang_rates}.
\begin{table}
\hspace*{-2cm}
\begin{tabular}{lll} \hline
Event & State transition & Transition rate \\ \hline
Infection of susceptible & $I_1 \to I_1 + 1$ & $(\beta/N) \left( \sum_{m=1}^k I_m \right)  \left( N - \sum_{m=1}^k I_m \right)$ \\
Transition to next infectious stage & $\left( I_{m-1} , I_m \right) \to \left( I_{m-1} - 1 , I_m + 1 \right)$ & $k \gamma I_{m-1} \mbox{ for } m=2,3,\ldots,k$ \\
Recovery of infected & $I_k \to I_k-1$ & $k \gamma I_k$ \\ \hline
\end{tabular}
\caption{Transition rates for the SIS model with $k$ infectious stages (Erlang distributed infectious periods).}\label{SIS_Erlang_rates}
\end{table}
The state space is $\left\{ \left( i_1 , i_2 , \ldots , i_k \right) \in \left( \Z_+ \right)^k : i_1 + i_2 + \cdots + i_k \le N \right\}$, and consists of a single absorbing state at $\bI = \left( I_1 , I_2 ,\ldots , I_k \right) = {\bf 0}$ together with a transient communicating class~$C$.
The basic reproduction number is $R_0 = \beta / \gamma$ as before.
The mean infectious period is $1/\gamma$ as before, but the variance of the infectious period is now $1 / k \gamma^2$ and thus may be varied according to our choice of~$k$.

We will consider only the case $R_0 > 1$, so that infection may become endemic in the population, and study the time to extinction from quasi-stationarity.
Denoting by $\bm{q}^k$ the quasi-stationary distribution, then the mean extinction time is given by $\tau_{\bq^k} = 1 / \lambda_k$ where $-\lambda_k$ is the eigenvalue with largest real part of~$Q_k^C$, the transition rate matrix restricted to the non-absorbing states.
The state-space is now of size ${N+k \choose k}$, so that finding the mean time to extinction from equation~(\ref{exact}) amounts to evaluating an eigenvalue of a sparse non-symmetric square matrix of dimension ${N+k \choose k} - 1$.

The process $\bI(t)/N$ may be approximated in the large-population limit by the deterministic process $\by = \left( y_1 , y_2 , \ldots , y_k \right)$ satisfying
\begin{eqnarray*}
{dy_1 \over dt} &=& \beta \left( 1 - \sum_{m=1}^k y_m \right) \left( \sum_{m=1}^k y_m \right) - k \gamma y_1 , \\
{dy_m \over dt} &=& k \gamma y_{m-1} - k \gamma y_m \mbox{ for } m=2,3,\ldots,k .
\end{eqnarray*}
This system has two equilibria: the disease-free equilibrium $\by = {\bf 0}$ and an endemic equilibrium $\by^* = \left( 1 - ( 1 / R_0 ) \right) ({\bf 1} / k)$, where $\bf 1$ denotes the vector with all entries equal to~1.

The (second-order) diffusion approximation is the process $\bY (t)$ satisfying
\begin{eqnarray*}
dY_1 &=& \left( {\beta \over N} \left( \sum_{m=1}^k Y_m \right) \left( N - \sum_{m=1}^k Y_m \right) - k \gamma Y_1 \right) dt\\
&& {}
+ \sqrt{{\beta \over N} \left( N - \sum_{m=1}^k Y_m \right) \left( \sum_{m=1}^k Y_m \right)} dW_0 - \sqrt{k \gamma Y_1} \, dW_1 ,\\
dY_m &=& k \gamma \left( Y_{m-1} - Y_m \right) dt
+ \sqrt{k \gamma Y_{m-1}} dW_{m-1} - \sqrt{k \gamma Y_m} \, dW_m \mbox{ for } m=2,3,\ldots,k,
\end{eqnarray*} 
where $W_0 , W_1 , W_2 , \ldots , W_k$ are independent standard Brownian motions.

Defining
\begin{eqnarray*}
\tau^{\mbox{Diff}} ( {\by } ) &=& E \left[ \left. \inf \left\{ t \ge 0 : \max_{1 \le m \le k} \left\{ Y_m(t) \right\} \le 0.5 \right\} \right| \bY (0) = \by \right]
\end{eqnarray*}
then $\tau^{\mbox{Diff}} ( \by )$ satisfies the Kolmogorov backward equation
\begin{eqnarray}
\sum_i A_i ( \bm{y} ) {\partial \tau^{\mbox{Diff}} ( \bm{y} ) \over \partial y_i} + {1 \over 2} \sum_{i,j} B_{ij} ( \bm{y} ) {\partial^2 \tau^{\mbox{Diff}} \over \partial y_i \, \partial y_j} &=& -1 \mbox{ in } \Omega , \label{KBE}
\end{eqnarray}
where the state-space is $\Omega = \left\{ \by \in ( \R_+ )^k : \max \left\{ y_1 , y_2 , \ldots , y_k \right\} \ge 0.5 , \; y_1 + y_2 + \cdots + y_k \le N \right\}$, and the coefficients are
\begin{eqnarray*}
\renewcommand{\arraystretch}{1.5}
\begin{array}{rcrcll}
&& A_1 ( \bm{y} ) &=& 
 {\beta \over N} \left( \sum_{m=1}^k y_m \right) \left( N - \sum_{m=1}^k y_m \right) - k \gamma y_1 ,\\
&& A_i ( \bm{y} ) &=& k \gamma \left( y_{i-1} - y_i \right) & i=2,3,\ldots,k ,\\
&& B_{11} ( \bm{y} ) &=& {\beta \over N} \left( \sum_{m=1}^k y_m \right) \left( N - \sum_{m=1}^k y_m \right) + k \gamma y_1 ,\\
&& B_{ii} (\bm{y})&=& k \gamma \left( y_{i-1} + y_i \right) & i=2,3,\ldots,k ,\\
B_{i,i-1} (\bm{y}) &=& B_{i-1,i} (\bm{y}) &=& - k \gamma y_{i-1} & i=2,3,\ldots,k ,\\
&& B_{ij} &=& 0 & j \notin \left\{ i-1,i,i+1 \right\} .
\end{array}
\renewcommand{\arraystretch}{1}
\end{eqnarray*}
The process has an absorbing boundary at $S_a = \left\{ \by \in \Omega : \max \left\{ y_1 , y_2 , \ldots , y_k \right\}  = 0.5 \right\}$, the remainder of the boundary, denoted $S_r$, being taken to be reflecting.
Denoting by $\bn = \left( n_1 , n_2 , \ldots , n_k \right)$ a vector normal to the boundary $S_r$, the boundary conditions (Gardiner~\cite{G09} section 6.6) are 
\begin{eqnarray*}
\tau^{\mbox{Diff}} ( \bm{y} ) &=& 0 \mbox{ on } S_a , \\
\sum_{i,j} n_i B_{ij} ( \bm{y} ) {\partial \tau^{\mbox{Diff}} ( \bm{y} ) \over \partial y_j} &=& 0 \mbox{ on } S_r .
\end{eqnarray*}

In this case, the boundary conditions on $S_r$ reduce to 
\begin{eqnarray*}
\renewcommand{\arraystretch}{3}
\begin{array}{rcll}
\displaystyle {\partial \tau^{\mbox{Diff}} \over \partial y_1} &=& 0 & \mbox{for } \by \in S_r \mbox{ with } y_1 = 0 , \\
\displaystyle {\partial \tau^{\mbox{Diff}} \over \partial y_i} 
- {\partial \tau^{\mbox{Diff}} \over \partial y_{i-1}} &=& 0 & \mbox{for } \by \in S_r \mbox{ with } y_i = 0 \quad (i=2,3,\ldots,k) , \\
\displaystyle {\partial \tau^{\mbox{Diff}} \over \partial y_k} &=& 0 & \mbox{for } \by \in S_r \mbox{ with } y_1 + y_2 + \cdots + y_k = N .
\end{array}
\renewcommand{\arraystretch}{1}
\end{eqnarray*}

The above system is not uniformly elliptic on $\Omega$, since the matrix $B ( \by )$ with entries $B_{ij} ( \by )$ has an eigenvalue of zero at points on the boundary.
Consequently, we have no rigorous proof that a unique solution exists.
Nevertheless, since our main interest is in approximations that work in practice, we will go ahead with investigating the performance of this diffusion approximation numerically.
Specifically, we apply the finite element method (see, for example,~\cite{R06}) using Freefem++ software~\cite{F++}.

In order to apply the finite element method, the partial differential equation must first be put into variational form, as follows.
We start from the integral formulation
\begin{eqnarray*}
\int_\Omega \left( \sum_i A_i ( \bm{y} ) {\partial \tau^{\mbox{Diff}} ( \bm{y} ) \over \partial y_i} + {1 \over 2} \sum_{i,j} B_{ij} ( \bm{y} ) {\partial^2 \tau^{\mbox{Diff}} \over \partial y_i \, \partial y_j} \right) w ( \bm{y} ) \, d\Omega &=& - \int_\Omega w ( \bm{y} ) \, d\Omega
\end{eqnarray*}
for appropriate test functions $w(\bm{y})$.
Integrating by parts in order to eliminate second order terms (treating mixed derivatives ${\partial^2 \over \partial y_i \partial y_j}$ symmetrically), applying boundary conditions to the resulting boundary integral terms, and noting that the condition $\tau^{\mbox{Diff}} ( \bm{y} ) = 0$ on $S_a$ implies that relevant test functions can be taken to satisfy $w ( \bm{y} ) = 0$ on $S_a$, then all boundary integral terms vanish, and we obtain the variational form
\begin{eqnarray}
{1 \over 2} \int_\Omega 
\sum_{i , j} \left( B_{ij} ( \bm{y} )  {\partial \tau^{\mbox{Diff}} \over \partial y_i} {\partial w \over \partial y_j}  
+ w(\bm{y}) {\partial B_{ij} \over \partial y_j} {\partial \tau^{\mbox{Diff}} \over \partial y_i} \right) d\Omega && \nonumber \\
{} - \int_{\Omega} \sum_i A_i ( \bm{y} ) w ( \bm{y} ) {\partial \tau^{\mbox{Diff}} \over \partial y_i} \, d\Omega 
&=& \int_\Omega w ( \bm{y} ) \, d\Omega .
\label{variational_form}
\end{eqnarray}
Freefem++ code to solve the above system in the case $k=2$ is given in the Appendix.

Approximating the drift and local variance matrices of the diffusion approximation to leading order around $N \by^*$ yields a multivariate Ornstein-Uhlenbeck approximation $\tilde {\bY}(t)$, centred at $N \by^*$, with 
drift matrix~$J$ and local variance matrix~$\tilde B$ given by
\begin{eqnarray*}
J &=& \left( \begin{array}{ccccc}
 2 \gamma - \beta - k \gamma & 2 \gamma - \beta &
 2 \gamma - \beta &
\cdots &  2 \gamma - \beta \\
k \gamma & - k \gamma & 0 & \cdots & 0 \\
0 & k \gamma & - k \gamma & \cdots & 0 \\
\vdots & \vdots & \vdots & & \vdots \\
0 & 0 & 0 & \cdots & -k \gamma
\end{array} \right) , \\
\tilde B &=& N \gamma \left( 1 - {1 \over R_0}\right) \left( \begin{array}{cccccc}
2  & - 1 & 0 & 0 & \cdots & 0 \\
- 1 & 2 & - 1 & 0 & \cdots & 0 \\
0 & - 1 & 2 & - 1 & \cdots & 0 \\
\vdots & \vdots & \vdots & \vdots && \vdots \\
0 & 0 & 0 & 0 & \cdots & 2 
\end{array} \right) .
\end{eqnarray*}
The stationary distribution of this Ornstein-Uhlenbeck process is multivariate Gaussian with mean $N {\by}^*$, variance matrix~$V$ satisfying
\begin{eqnarray}
J V + V^T J &=& - \tilde B . \label{Sigma}
\end{eqnarray}
The solution to equation~(\ref{Sigma}) is
\begin{eqnarray*}
V &=& {N \over k} \left( 1 - {1 \over R_0} \right) I_k + {N \over k^2}  \left( {2 \over R_0} - 1 \right) \mathbbm{1}_k ,
\end{eqnarray*}
where $I_k$ is the $k \times k$ identity matrix and $\mathbbm{1}_k$ represents the $k \times k$ matrix with all entries equal to~1.

Under the Ornstein-Uhlenbeck approximation, the total number of infectives $\tilde Y_1 + \tilde Y_2 + \cdots + \tilde Y_k$ has a Gaussian stationary distribution with mean $N y^*$ and variance $N / R_0$.
That is, the Ornstein-Uhlenbeck approximation suggests that the quasi-stationary distribution of the total number of infected individuals does not depend upon~$k$, the number of infected stages.
Further, under the Ornstein-Uhlenbeck approximation $\tilde Y_1 , \tilde Y_2 , \ldots , \tilde Y_k$ are exchangeable random variables, in stationarity.
Hence we can approximate
\begin{eqnarray*}
q^k_{\be_k} &\approx& {q_1 \over k} 
\end{eqnarray*} 
where $\be_k$ denotes the unit vector with $k$th element equal to~1, $q^k_{\be_k}$ is the quasi-stationary probability of being in state~$\be_k$, and $q_1$ is the quasi-stationary probability for the classic ($k=1$) SIS model to be in state $I=1$.
The mean time to extinction $\tau_{\bq^k}$ is then approximated as
\begin{eqnarray*}
\tau_{\bq^k} &=& {1 \over k \gamma q^k_{\be_k}} 
\approx {1 \over \gamma q_{1}} 
\approx \tau^{\mbox{OU}} ,
\end{eqnarray*}
where $\tau^{\mbox{OU}}$ is given by equation~(\ref{OU}), exactly as for the classic SIS model with exponentially distributed infectious periods.

Moving on to the Hamiltonian approach, our SIS model with Erlang distributed infectious periods has Hamiltonian 
\begin{eqnarray*}
H_k ( \by , \btheta ) &=& \beta \left( \sum_{m=1}^k y_m \right) \left( 1 - \sum_{m=1}^k y_m \right) \left( {\rm e}^{\theta_1} - 1 \right) \\
&& {} + k \gamma \sum_{m=1}^{k-1} y_{m} \left( {\rm e}^{-\theta_{m} + \theta_{m+1}} - 1 \right)
+ k \gamma y_k \left( {\rm e}^{-\theta_k} - 1 \right) ,
\end{eqnarray*}
describing a motion in $\R^{2k}$ that follows the equations of motion
\begin{eqnarray*}
{dy_i \over dt} &=& {\partial H_k \over \partial \theta_i} \mbox{ for } i=1,2,\ldots,k,\\
{d\theta_i \over dt} &=& - {\partial H_k \over \partial y_i} \mbox{ for } i=1,2,\ldots,k.
\end{eqnarray*}

Setting $\theta_0 = \theta_{k+1} = y_0 = 0$, these equations of motion may be written as
\begin{eqnarray}
\left.
\begin{array}{rcl}
\displaystyle {dy_i \over dt} &=& 
\beta \left( \sum_{m=1}^k y_m \right) \left( 1 - \sum_{m=1}^k y_m \right) {\rm e}^{\theta_1} \delta_{i1} + k \gamma \left( y_{i-1} {\rm e}^{-\theta_{i-1}+\theta_i} - y_i {\rm e}^{-\theta_i+\theta_{i+1}} \right)  
\\
&& \hspace{8cm} \mbox{ for } i=1,2,\ldots,k,\\
\displaystyle {d\theta_i \over dt} &=& 
- \beta \left( 1 - 2 \sum_{m=1}^k y_m \right) \left( {\rm e}^{\theta_1} - 1 \right)
- k \gamma \left( {\rm e}^{-\theta_i + \theta_{i+1}} - 1 \right)
\mbox{ for } i=1,2,\ldots,k,
\end{array} \right\}
\label{Hamiltonian_system_k}
\end{eqnarray}
where $\delta_{ij}$ is the Kronecker delta.
This system has classical equilibrium points at $\left( \bm{y} , \bm{\theta} \right) = \left( {\bf 0} , \bm{0} \right)$ and $\left( \by^* , \bm{0} \right)$, together with a non-classical disease-free equilibrium point at $\left( \bm{0} ,   \bm{\theta}^* \right)$, where $\bm{\theta}^* = ( k, k-1, \ldots , 3,2,1) \theta_k^*$ with $\theta_k^*$ satisfying
\begin{eqnarray}
{\beta \over k \gamma} \left( \sum_{m=1}^k {\rm e}^{m \theta_k^*} \right) - 1 &=& 0 . \label{theta_k}
\end{eqnarray}
The function $f(z) = \sum_{m=1}^k z^m$ is increasing, with $f(0) = 0$ and $f(z) \to \infty$ as $z \to \infty$, so there is a unique positive $z^*$ satisfying $f(z^*) = k \gamma / \beta$, and then the unique real solution of equation~(\ref{theta_k}) is $\theta_k^* = \ln z^*$.
Also, for $R_0 > 1$ we have $f(z^*) < k$, and  since $f(1) = k$ it follows that $z^* < 1$, so $\theta_k^* < 0$.

We may approximate $\ln \left( \tau_{\bq} \right) / N$ by $\ln \left( \tau_{H_k} \right) / N$, where
\begin{eqnarray}
\tau_{H_k} &=& \exp \left( N A_k \right) \label{tau_H_k}
\end{eqnarray}
with action~$A_k$ given by
\begin{eqnarray}
A_k &=& \int_{-\infty}^\infty \sum_{m=1}^k \theta_k  {dy_k \over dt} \, dt , \label{A_k}
\end{eqnarray}
the integration being along a trajectory going from $\left( \by^* ,\, {\bf 0} \right)$ at $t=-\infty$ to $\left( {\bf 0} , \, \bm{\theta}^* \right)$ at $t=+\infty$.

Solving the $2k$-dimensional system of ordinary differential equations~(\ref{Hamiltonian_system_k}) numerically for the case $k=2$, subject to boundary conditions at $t=-\infty$ and $t=+\infty$, and then using the solution $\left( \by(t) , \bm{\theta} (t) \right)$ thus obtained to evaluate $A_k$, it appears that $A_2 = A = (1/R_0) - 1 + \ln R_0$.
Now formula~(\ref{A_k}) is equivalent to $A_k = - S_k ( \btheta^* )$, where the function $S_k ( \btheta )$ satisfies the partial differential equation $H_k \left( {\partial S_k \over \partial \btheta} , \btheta \right) = 0$ with $S_k ( \bm{0} ) = 0$.
From the form of the equilibrium point $\btheta^*$, together with the conjecture that $A_k = A$, it is possible to guess the solution to be
\begin{eqnarray}
S_k ( \bm{\theta} ) &=& \ln \left( {\sum_{m=1}^k {\rm e}^{\theta_m} \over k} \right) - {\gamma \over \beta} \left( 1 - {k  \over \sum_{m=1}^k {\rm e}^{\theta_m}} \right) . \label{S_k}
\end{eqnarray}
It is now straightforward to check that $S_k ( \bm{\theta} )$ given by formula~(\ref{S_k}) does indeed satisfy $H_k \left( {\partial S_k \over \partial \btheta} , \btheta \right) = 0$, and that $S_k ( \btheta^* ) = 1 - (1/R_0) - \ln R_0$.
Hence $A_k = (1/R_0) - 1 + \ln R_0$ for $k=1,2,\ldots$.
We thus arrive at the same conclusion as suggested by the Ornstein-Uhlenbeck approximation: that, to leading order, the mean persistence time does not depend upon~ $k$.

For the SIS model with Erlang distributed infectious periods, a precise result corresponding to formula~(\ref{AD}) is available from the recent paper~\cite{BBN16}.
Consider an SIS epidemic model in which individual infectious periods are distributed as any non-negative random variable~$Q$ of finite variance, and suppose that $R_0 = \beta E[Q] > 1$.
Denote by $p_Q$ the asymptotic (large $N$) probability, starting from a single infected individual in an otherwise susceptible population, that only a minor outbreak occurs.
Denoting $g( \theta ) = E \left[ {\rm e}^{\theta Q} \right]$ then it is well-known that $p_Q$ is the unique solution in $[0,1)$ of
\begin{eqnarray*}
p_Q &=& g \left( - \beta \left( 1 - p_Q \right) \right) .
\end{eqnarray*}
Lemmas~3.2 and~3.3 of~\cite{BBN16} together imply that $\tau \sim \tau^{\mbox{BBN}}$, where
\begin{eqnarray}
\tau^{\mbox{BBN}} &=& E[Q] \, \sqrt{{2\pi \over N}} {1 \over \left( R_0 - 1 \right) \left( 1 - p_Q \right)} \, \exp \left( N \left( (1/R_0) - 1 + \ln R_0 \right) \right) . \label{BBN}
\end{eqnarray}
In the case of the classic SIS model, then $p_Q = 1 / R_0$ and formula~(\ref{AD}) is recovered.
When $Q$ follows an Erlang distribution with mean $E[Q] = 1/\gamma$ and variance $1/k\gamma^2$, then $g(\theta) = \left( 1 - ( \theta / k \gamma ) \right)^{-k}$, so that $p_Q$ is the unique solution in $[0,1)$ of
\begin{eqnarray}
p_Q \left( 1 + R_0 ( 1- p_Q ) / k \right)^k &=& 1 . \label{pQ}
\end{eqnarray}

\subsection{Numerical results}
For numerical work, we take $k=2$, and investigate how well our general approximation methods perform for a 2-dimensional problem.

\begin{center}
\begin{figure}
\includegraphics[width=12cm]{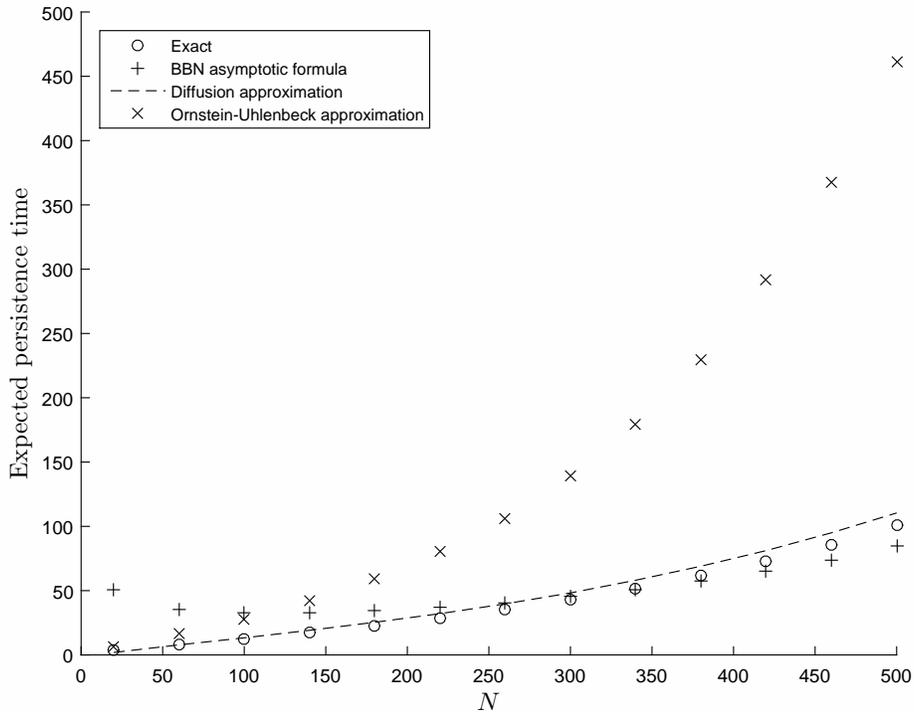}
\caption{Expected persistence time from quasi-stationarity ($\tau^{\bm{q}}$) for the SIS model with Erlang distributed infectious periods having shape parameter $k=2$, in a population of size $N$, together with approximations.
Parameter values $\beta = 1.1$, $\gamma = 1$ (so that $R_0 = 1.1$ and $y_1^* = y_2^* \approx 0.0455$).
Exact values computed from formula~(\ref{QSD}); diffusion approximation~(\ref{KBE}); Ornstein-Uhlenbeck approximation~(\ref{OU}).
Diffusion approximation is initiated at the finite element meshpoint nearest to $\bY(0) = N \left(  y_1^* , y_2^* \right)$.}
\label{fig12}
\end{figure}
\end{center}

Consider first a system not far above threshold, with $R_0 = 1.1$.
Comparing figure~\ref{fig12} with figure~\ref{fig5}, the pictures seem rather similar, but there are significant differences in how they are constructed.
In both cases, the exact values are computed by solving equation~(\ref{QSD}).
For the classic SIS model this requires the evaluation of an eigenvalue of a square matrix $Q_C$ of dimension~$N$; for the case $k=2$, the matrix in question is of dimension~$N(N + 3)/2$.

The diffusion approximation for the classic SIS model is computed by evaluating the integral~(\ref{Diff}).
For $k=2$, no such explicit expression for $\tau^{\mbox{Diff}}$ is available, and we  proceed via the Finite Element Method.
Freefem++ code (see Appendix) is used, which gives a numerical approximation to $\tau^{\mbox{Diff}} ( \by )$ for a grid of $\by$ values within $\Omega$.
The user specifies the set of grid points along the boundary $S_a \cup S_r$, and then grid points in the interior of $\Omega$ are selected automatically by the software.
To approximate $\tau_{\bm{q}^k}$, we select the grid point at minimal Euclidean distance from $N \by^*$.
Since the surface $\tau^{\mbox{Diff}} ( \by )$ is quite flat around $N \by^*$, our approximation is not overly sensitive to the precise value of $\by$ chosen.

Turning to the asymptotic result of~\cite{BBN16}, with $k=2$ then $p_Q$ satisfying equation~(\ref{pQ}) may be found explicitly as
\begin{eqnarray*}
p_Q &=& {4 + R_0 - \sqrt{R_0 \left( 8 + R_0 \right)} \over 2 R_0} .
\end{eqnarray*}
Formula~(\ref{BBN}) in this case reduces to
\begin{eqnarray}
\tau^{\mbox{BBN}} &=& {1 \over \gamma} \, \sqrt{{2\pi \over N}} {2 R_0 \exp \left( N \left( (1/R_0) - 1 + \ln R_0 \right) \right)\over \left( R_0 - 1 \right) \left( R_0 - 4 + \sqrt{R_0 \left( 8 + R_0 \right)} \right)}  . \label{BBN}
\end{eqnarray}

\begin{center}
\begin{figure}
\includegraphics[width=12cm]{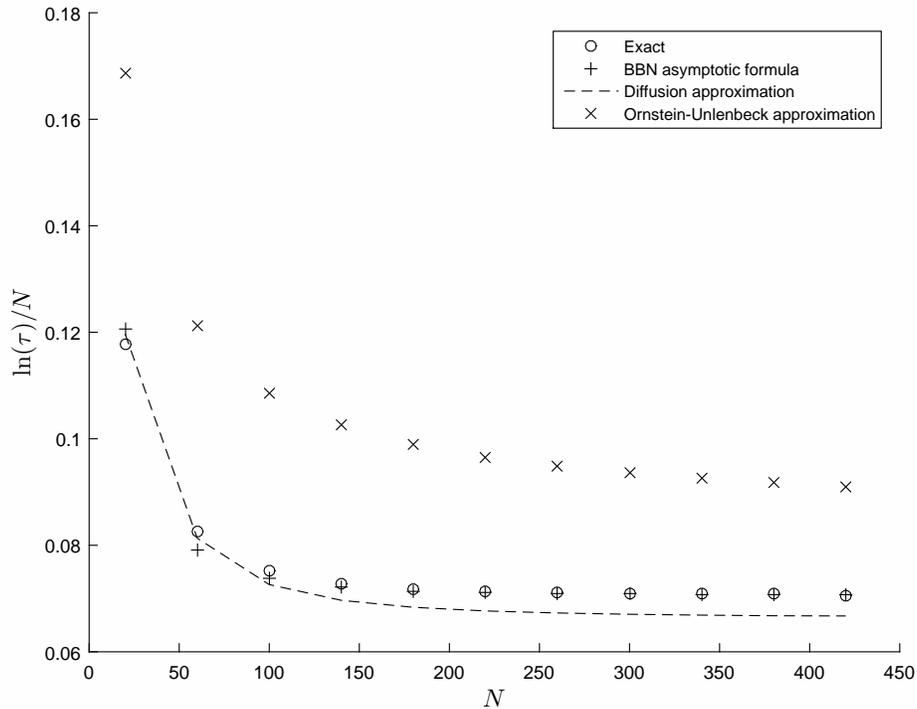}
\caption{$\ln ( \tau_{\bm{q}} ) / N$ for the SIS model with Erlang distributed infectious periods having shape parameter $k=2$ in a population of size $N$, together with approximations.
Parameter values $\beta = 1.5$, $\gamma = 1$ (so that $R_0 = 1.5$ and $\by^* \approx ( 0.1667 , 0.1667 )$).
Exact values computed from formula~(\ref{QSD}); diffusion approximation~(\ref{KBE}); Ornstein-Uhlenbeck approximation~(\ref{OU}).
Diffusion approximation is initiated at the finite element meshpoint nearest to $\bY(0) = N \left(  y_1^* , y_2^* \right)$.}
\label{fig16}
\end{figure}
\end{center}

The performance of our approximations for a system well above threshold ($R_0 = 1.5$) is shown in figure~\ref{fig16}.
The picture here is very similar to that seen in figure~\ref{fig6} for the classic SIS model: $\tau^{\mbox{BBN}}$ gives the correct leading-order asymptotic behaviour; $\tau^{\mbox{Diff}}$ underestimates mean persistence time; and $\tau^{\mbox{OU}}$ provides a very poor approximation.
Correct leading-order asymptotic behaviour, in agreement with the result of~\cite{BBN16}, may be found using the Hamiltonian approach.

\section{Discussion}
For the below-threshold case $R_0 < 1$, we have seen that a linear birth-death process approximation works well provided the initial number of infected individuals is small.
Although we only looked at the classic SIS model, a wide range of infection models \cite{B14} are known to satisfy an epidemic threshold theorem.
The threshold theorem tells us that for $R_0 < 1$ and population size $N$ large, infection is very likely to die away quickly without ever infecting a significant fraction of the population; and that throughout a typical short-lived outbreak the number of susceptible individuals $S(t)$ may be well approximated by the constant $N$, so that consequently the process of infected individuals $I(t)$ may be well-approximated by a (linear) branching process.
Thus we have a quite general approximation technique, which is seen to perform well in practice provided the number of initially infected individuals is small (figures~\ref{fig1},~\ref{fig2}).
A number of previous authors have also considered the case in which a significant fraction of the population is initially infected.
This situation is of less practical interest, since for sufficiently large~$N$ the chance of a significant fraction of the population ever becoming infected is negligible.
The only approximation that we have found to work well in this case is the `diffusion approximation' (figure~\ref{fig4}), and the value of this approximation seems doubtful, since it appears no more straightforward to compute than the exact answer (indeed, if anything it seems more difficult to work with).

We turn now to the more interesting case $R_0 > 1$, so that it is possible for infection to sustain itself in the population over the long term.
In this case, we studied time to extinction from an endemic state.
First of all, consider the `diffusion approximation'.
We saw that for moderate population sizes and $R_0$ not too far above~1, the approximation performs rather well (figures~\ref{fig5},~\ref{fig12}).
On the other hand, it is known that this approximation gives the wrong asymptotic (large~$N$) behaviour for the classic SIS model~\cite{DSS05}.
We have seen (figures~\ref{fig6},~\ref{fig16}) that for the SIS model with Erlang-distributed infectious periods, exactly as for the classic SIS model, the diffusion approximation gives the wrong asymptotic behaviour.
Hence we conclude that, firstly, the diffusion approximation does not seem trustworthy in cases where it is not possible to compute the exact answer for comparison; and secondly, the circumstances in which the diffusion approximation performs well seem to be precisely those circumstances in which it is straightforward to numerically evaluate the exact solution~(\ref{exact}).
Thus the value of this approximation method seems at best doubtful.
The recent paper~\cite{WGPAI14} uses a diffusion approximation to study time to extinction of E. coli O157:H7 infection from a herd of cattle.
The diffusion is approximating a 2-dimensional Markov jump process, with population size of the order of $100$~cattle and $10^{11}$~colony forming units (cfu) of E. coli.
The state-space of the Markov jump process is therefore far too large to compute the mean extinction time exactly from equation~(\ref{exact}), and some approximation technique is required.
For different scenarios considered in~\cite{WGPAI14}, mean extinction times computed from the diffusion approximation
are from around 1~month up to around 40~years.
This suggests that the infection process is not very far above threshold, and that the approximation may perform well.
On the other hand, a population size of $10^{11}$~cfu is very large, which may lead the approximation to break down.
Although it is certainly possible in these circumstances that the diffusion may give a good approximation to mean extinction time, our results suggest that one cannot have any real confidence of this.
Similarly, in~\cite{DV16}, diffusion processes are used to compute mean extinction times for a Lotka-Volterra predator-prey model, for two variants of the SIS epidemic model, and for a two-pathogen epidemic model; there is no attempt to validate the approximations with reference to the underlying Markov chains.

Secondly, consider the Ornstein-Uhlenbeck approximation.
This approach has been widely used in the literature due to its great simplicity. 
A formula such as~(\ref{OU}) is much more analytically tractable than the exact solution~(\ref{Norden}) or the diffusion approximation~(\ref{Diff}).
However, we have seen (figures~\ref{fig5},~\ref{fig6},~\ref{fig12},~\ref{fig16}) that this approximation performs extremely poorly.
Nevertheless, we would argue that this approach retains some value in limited circumstances: specifically, for qualitative comparison between two infection models which share the same deterministic endemic prevalence level.
This is because, for $N$~large and $R_0$ well above~1, this approach gives a good approximation to the body of the quasi-stationary distribution --- see, for example,~\cite{C05,N99}.
The approximation fails in the tails of the distribution, hence the poor performance of our numerical approximation, which relies upon approximating the quasi-stationary probability that exactly one individual is infected.
For qualitative comparison of two models, provided that both models have the same mean number of infected individuals in quasi-stationarity (approximated by the deterministic endemic level), it seems reasonable to suppose that that the model whose quasi-stationary distribution has higher variance will assign greater probability to the state with one individual infected.
The Ornstein-Uhlenbeck approach gives a good approximation for the variance of the quasi-stationary distribution, provided $N$ is large and $R_0 \gg 1$.
It is therefore reasonable to suppose that for two infection models sharing a common deterministic endemic prevalence level, the model under which the variance of the stationary distribution of the Ornstein-Uhlenbeck approximating process is greater will have the lower mean persistence time.

Finally, we turn to the Hamiltonian approach.
Whereas diffusion approximations deal with moderate deviations from a deterministic mean, the Hamiltonian approach is suited to dealing with large deviations, and so would be expected to perform better when studying time to disease extinction.
Our results (figures~\ref{fig6},~\ref{fig16}) confirm that this approach does indeed give the correct leading-order asymptotic behaviour for large~$N$.
This approach therefore seems the most promising overall, but a number of difficulties remain.
Firstly, we have not obtained a useful numerical approximation for finite~$N$, due to the unknown pre-factor $C$ in the relation $\tau_{\bq} \sim C ( N , R_0 ) \exp (NA)$.
For 1-dimensional systems, including the classic SIS model, it is possible to evaluate the pre-factor~$C$ by retaining higher order information when approximating $K ( \theta , t )$, see~\cite{DSS05,AM10}.
It does not seem straightforward to generalise this to higher dimensional systems, although some progress has been made for one particular 2-dimensional infection model in~\cite{VG95}, where is is shown that the pre-factor is of the form $C_0 / \sqrt{N}$ for some (unknown) $C_0$ that does not depend upon~$N$.
Secondly, to evaluate the constant $A$ in the dominant exponential term will, in general, require the solution of a boundary value problem for a $2k$ dimensional ordinary differential equation system, where $k$ is the dimensionality of the original model.
Recent work on numerical approaches to this problem includes~\cite{SFBS11,LSS14}, where systems in $k=2, 3$ dimensions are analysed.
An alternative is to seek approximations valid within certain regions of parameter space, such as the `adiabatic approximation' of~\cite{DSL08}, where (for their model of interest) a simple explicit expression for $A$ valid for $R_0$ close to~1 is derived.
Our results for the SIS model with Erlang-distributed infectious periods suggest another possibility: that for systems with sufficient symmetry, it may be possible to guess an explicit solution to the partial differential equation satisfied by $S ( \bm{\theta} )$, as we did to obtain formula~(\ref{S_k}), which leads directly to an explicit formula for $A$.
In fact, it is not necessary to solve for $S ( \bm{\theta} )$ for all $\bm{\theta}$, but only to evaluate $S \left( \bm{\theta}^* \right) - S \left( {\bf 0} \right)$.
To illustrate this, consider an SEIS model, with a latent (`exposed') period between being infected and becoming infectious.
Suppose the latent period follows an Erlang distribution corresponding to $j$ latent stages, each of mean $1 / j \nu$, while the infectious period remains Erlang with $k$ stages each of mean $1 / k\gamma$.
The Hamiltonian for this system is
\begin{eqnarray*}
H_{j,k} ( \by , \btheta ) &=& \beta \left( \sum_{m=j+1}^{j+k} y_m \right) \left( 1 - \sum_{m=1}^{j+k} y_m \right) \left( {\rm e}^{\theta_1} - 1 \right) \\
&& {} + j \nu \sum_{m=1}^j y_m \left( {\rm e}^{-\theta_m + \theta_{m+1}} - 1 \right) \\
&& {} + k \gamma \sum_{m=j+1}^{j+k-1} y_{m} \left( {\rm e}^{-\theta_{m} + \theta_{m+1}} - 1 \right)
+ k \gamma y_{j+k} \left( {\rm e}^{-\theta_{j+k}} - 1 \right) 
\end{eqnarray*}
where $y_1 , y_2 , \ldots , y_j$ correspond to latent stages and $y_{j+1} , y_{j+2} , \ldots , y_{j+k}$ to infectious stages, and similarly for the components of $\bm{\theta}$.
The corresponding equations of motion have a disease-free equilibrium point with $\bm{\theta}^* = ( k,k,\ldots,k,k-1,k-2,\ldots,3,2,1) \theta_k^*$, where $\theta_k^*$ is given by equation~(\ref{theta_k}).
On the hyperplane $\theta_1 = \theta_2 = \cdots = \theta_j = \theta_{j+1}$, the function $S_k ( \bm{\theta} )$ given by formula~(\ref{S_k}) (with components $\theta_1 , \theta_2 , \ldots , \theta_k$ re-labelled as $\theta_{j+1} , \theta_{j+2} , \ldots , \theta_{j+k}$) is readily seen to satisfy the relevant Hamilton-Jacobi equation, and it follows that the action is given by $A_{j,k} = (1/R_0) - 1 + \ln R_0$ as before.
Thus we see that the existence of a latent period has no effect upon the leading-order term in the expected time to extinction of infection.
The above argument can be straightforwardly extended to allow for a latent period distributed according to any phase-type distribution.

\section*{Acknowledgements}
Elliott Tjia was supported by a studentship from the Engineering and Physical Sciences Research Council.
The authors would like to thank Bernd Schroers, Robert Weston and Des Johnston for helpful discussions regarding the Hamiltonian approach.

\section*{Appendix}
Freefem++ code~\cite{F++} to numerically solve the Kolmogorov backward equation~(\ref{variational_form}) in the case $k=2$.

\begin{verbatim}
// set parameter values
real gamma=1., beta=2.0, n=260., epsilon=0.5; 

// define the boundary
border B1a(t=0,epsilon){x=t; y=epsilon;}
border B1b(t=0,epsilon){x=epsilon; y=epsilon-t;}
border B2(t=epsilon,n){x=t; y=0;}
border B3(t=0,n){x=n-t; y=t;}
border B4(t=0,n-epsilon){x=0; y=n-t;}

// define mesh
mesh Th = buildmesh (B1a(10)+B1b(10)+B2(100)+B3(100)+B4(100));

// define finite element space
fespace Vh(Th,P2);
Vh tau,w; 

// solve PDE
solve Backward(tau,w,solver=LU) = 
int2d(Th)(-((beta/(2*n))*(x+y)*(n-x-y)+gamma*x)*dx(tau)*dx(w)
-gamma*(x+y)*dy(tau)*dy(w)
+gamma*x*dy(tau)*dx(w)
+gamma*x*dx(tau)*dy(w))
+int2d(Th)(2*gamma*(x-y)*w*dy(tau)
+((beta/n)*(x+y)*(n-x-y)-2*gamma*x-(beta/(2*n))*(n-2*x-2*y)-gamma)*w*dx(tau))
-int2d(Th)(-w)
+ on(B1a,tau=0) + on(B1b,tau=0) ; 

// Output results, for processing in Matlab
{ ofstream ff("SISk.txt");
for (int i=0;i<Th.nt;i++)
{ for (int j=0; j <3; j++)
ff<<Th[i][j].x << " "<< Th[i][j].y<< " "<<tau[][Vh(i,j)]<<endl;
ff<<Th[i][0].x << " "<< Th[i][0].y<< " "<<tau[][Vh(i,0)]<<"\n\n\n";
}
}
\end{verbatim}

\end{document}